# Energy Efficiency Optimization of Multi-unit System with Different Devices

Fulai Yao
School of Artificial Intelligence
Hebei University of Technology
Tianjin, China
*fulaiyao@aliyun.com

**Abstract-** The energy efficiency optimization of the power generation system and the energy efficiency optimization of the energy consumption system are unified into the same optimization problem, and a simple method to achieve energy efficiency optimization without establishing an accurate mathematical model of the system is proposed. For systems with similar energy efficiency, it is proved that the best load distribution method between equipment is to keep the operating energy efficiency of each operating device equal, Yao's theorem 1. It is proved that the optimal switching method for the number of operating units between equipment with different energy efficiency is to keep the energy efficiency of the switching point equal, or at the maximum load point of the equipment, Yao's Theorem 2. This article gives two cases, a system composed of equipment with similar efficiency and a system composed of equipment with different efficiency.

## 1.1 Introduction

There are many power stations in the world with the same or different types of generators, and there are many energy-consuming systems with the same or different types of equipment.

As long as multiple devices are required to cooperate to complete a job, there must be an energy-saving problem. This is the problem of overall energy efficiency optimization.

In a power plant, for a given raw material input, there is a maximum energy output optimization problem where we should decide how many generators to use and how much load each generator should bear. For energy-consuming systems, for a given energy output, there is a minimum energy input optimization problem where we should decide how many devices to use and how much load each device should be assigned.

There are many optimization methods that have been widely studied, such as genetic algorithm [1], ant colony optimization [2], recursive quadratic programming [3], equal increment principle [4], linear programming [5], and stochastic dynamic programming [6], neural network [7], evolutionary strategy [8], Lagrangian relaxation [9], particle swarm optimization [10], simulated annealing [11], quasi-Newton [12], etc.[15]-[19] Many optimization methods require the establishment of system models. Since accurate models of many actual systems are difficult to establish, these methods become difficult to apply. For a system composed of equipment with the same efficiency, a method to obtain optimal operating energy efficiency without establishing an accurate mathematical model of the system was proposed [13]. For a system composed of equipment with similar efficiencies, an energy efficiency optimal control method is proposed [14].

This paper further proposes an energy efficiency optimization control method for systems composed of equipment with different efficiencies, which can avoid the difficulty of system modeling.

### 1.2 Efficiency Optimization of Multi-unit System with Different Devices

A multi-unit system is a system composed of multiple devices to accomplish the same task. For example, a power plant composed of multiple generators, a power distribution station composed of multiple transformers, a pumping station composed of multiple pumps, a high-speed rail train driven by multiple motors, and a hydrogen production station composed of multiple hydrogen generators.

For a fixed total energy-A input $P_t$, the maximum output energy B $W_t$ is expressed as:

$$W_t = \max \sum_{i=1}^{n} P_i \eta_{Pi}(P_i) \\ s.t.\ \sum_{i=1}^{n} P_i = P_t \\ P_{im} \geq P_i \geq 0 \tag{1}$$

For a fixed total energy-B output $W_t$, the minimization of the input energy A total $P_t$ can be expressed as:

$$P_t = min \sum_{i=1}^{n} \frac{W_i}{\eta_{Wi}(W_i)} \\ s.t.\ \sum_{i=1}^{n} W_i = W_t \\ W_{im} \geq W_i \geq 0 \tag{2}$$

The energy efficiency functions $\eta_{Pi}$ ($P_i$) and $\eta_{Wi}(W_i)$ are not the identical function. They have different variables and different function values.

Equations (1) and (2) are two different expressions of the same optimization problem, both of which are operational energy efficiency optimization problems. Just solve one of them and the problem is solved. We choose to solve the optimal result of equation (1).

It can be seen from equation (1) that since the energy efficiency function is nonlinear and has restrictions on the maximum input value $P_{im}$, this is a constrained nonlinear optimization problem. Since the number of running units n is an integer, the optimal load distribution value $P_i$ is a real number, which is another integer-real mixed optimization problem.

### 1.3 Energy Efficiency Function

The energy efficiency functions have a similar shape. Let's take $\eta_{Pi}(P)$ as an example. For convenience of writing, we use $\eta_i(P)$ instead of $\eta_{Pi}(P)$. The energy efficiency curve of $\eta_i(P)$ is shown in Fig. 1.

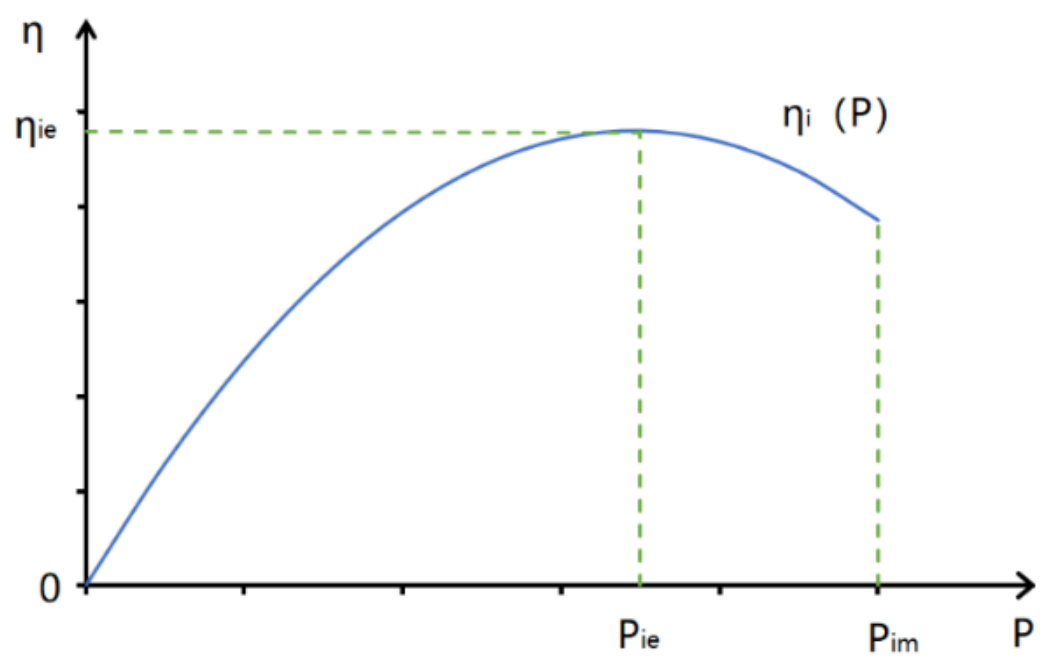

Fig. 1 Energy efficiency curve $\eta_i(P)$

In Fig. 1, $P_{im}$ is the maximum A energy input of the i-th device, $\eta_{ie}$ is the highest operating energy efficiency of the i-th device, $P_{ie}$ is the A energy input of the i-th device at the highest operating energy efficiency, and $\eta_i(P)$ has features:

$$0 \leq \mathrm{P} \leq P_{im}$$
$$\eta_{ie} = \eta_i(P_{ie})$$
$$0 \leq \eta_i(P) \leq \eta_{ie} \tag{3}$$
$$\eta_i(0) = 0$$
$$\eta_i''(P) < 0$$

Assume the efficiency curves $\eta_1(P)$ and $\eta_i(P)$ of the 1st and ith devices are shown in Fig. 2.

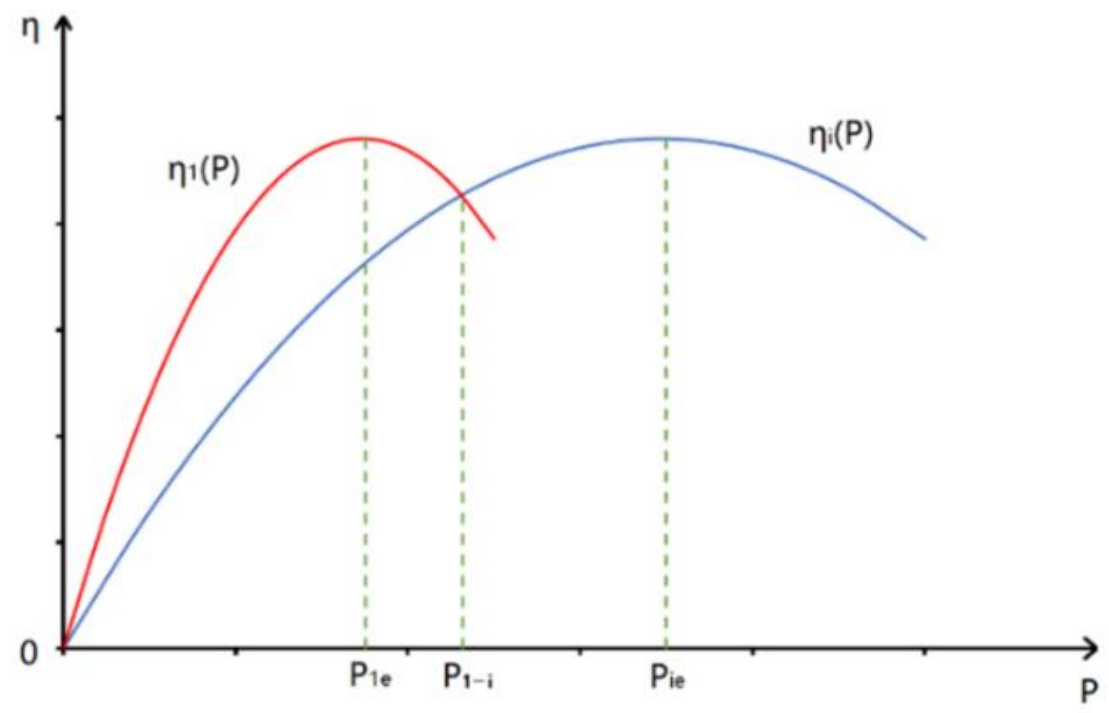

Fig. 2 Energy efficiency curves $\eta_1(P)$ and $\eta_i(P)$

If the following equation holds:

$$\eta_i(P_i) = \eta_1\left(\frac{P_1}{\beta_i}\right) \tag{4}$$

Among them, $\beta_i$ is a constant. We refer to the i-th device and the first device as devices with similar efficiency, called "similar efficiency devices".

i=1, then the i-th device and the first device are the device with the identical efficiency, recorded as $\beta_1=1$;

Equipment with similar efficiency has the following characteristics:

$$\eta_i'(P_i) = \frac{\eta_1'(P_1)}{\beta_i}$$
$$\beta_i = \frac{\eta_1'(P_1)}{\eta_i'(P_i)} \tag{5}$$
$$P_{ie} = \beta_i P_{1e}$$
$$\beta_i = \frac{P_{ie}}{P_{1e}}$$

**1.4 Optimal Load Distribution Theorem (Yao Theorem 1)**

A multi-unit system has n devices running, and each device has a different efficiency function. Assume $P_i$ is greater than 0, that is, every operating device works. Equation (1) simplifies to

$$W_t = \max \sum_{i=1}^{n} P_i \eta_i(P_i)$$
$$s.t.\ \sum_{i=1}^{n} P_i = P_t > 0 \tag{6}$$
$$P_{imax} \geq P_i > 0$$

Assuming n is optimal, consider the following three situations:

1) n=2

System C has two variables $P_1$ and $P_2$

$$P_1 + P_2 = P_t$$
$$P_1 > 0 \tag{7}$$
$$P_2 > 0$$

The objective function can be expressed as

$$W_t = P_1\eta_1(P_1) + P_2\eta_2(P_2) \tag{8}$$

The optimization condition is

$$W_t'(P_1) = 0 \tag{9}$$

According to known conditions

$$P_2 = P_t - P_1 \tag{10}$$

We have

$$\eta_1(P_1) + P_1\eta_1'(P_1) - \eta_2(P_2) - P_2\eta_2'(P_2) = 0 \tag{11}$$

a) If $\eta_1(P_1) = \eta_2(P_2)$, and satisfy the following equation

$$P_1\eta_1'(P_1) = P_2\eta_2'(P_2) \tag{12}$$

The two derivatives have the same signs and are both in the rising or falling side of the efficiency curve, then equation (11) is established. Equation (12) shows if abs ($\eta_1'(P_1)$) is large, $P_1$ is small, and if abs ($\eta_1'(P_1)$) is small, $P_1$ is large.

If the two devices are with similar efficiencies, $\eta_1(P_1) = \eta_2(P_2)$, $P_1\eta_1'(P_1) = P_2\eta_2'(P_2)$, and equation (11) is true. As shown in Fig. 3.

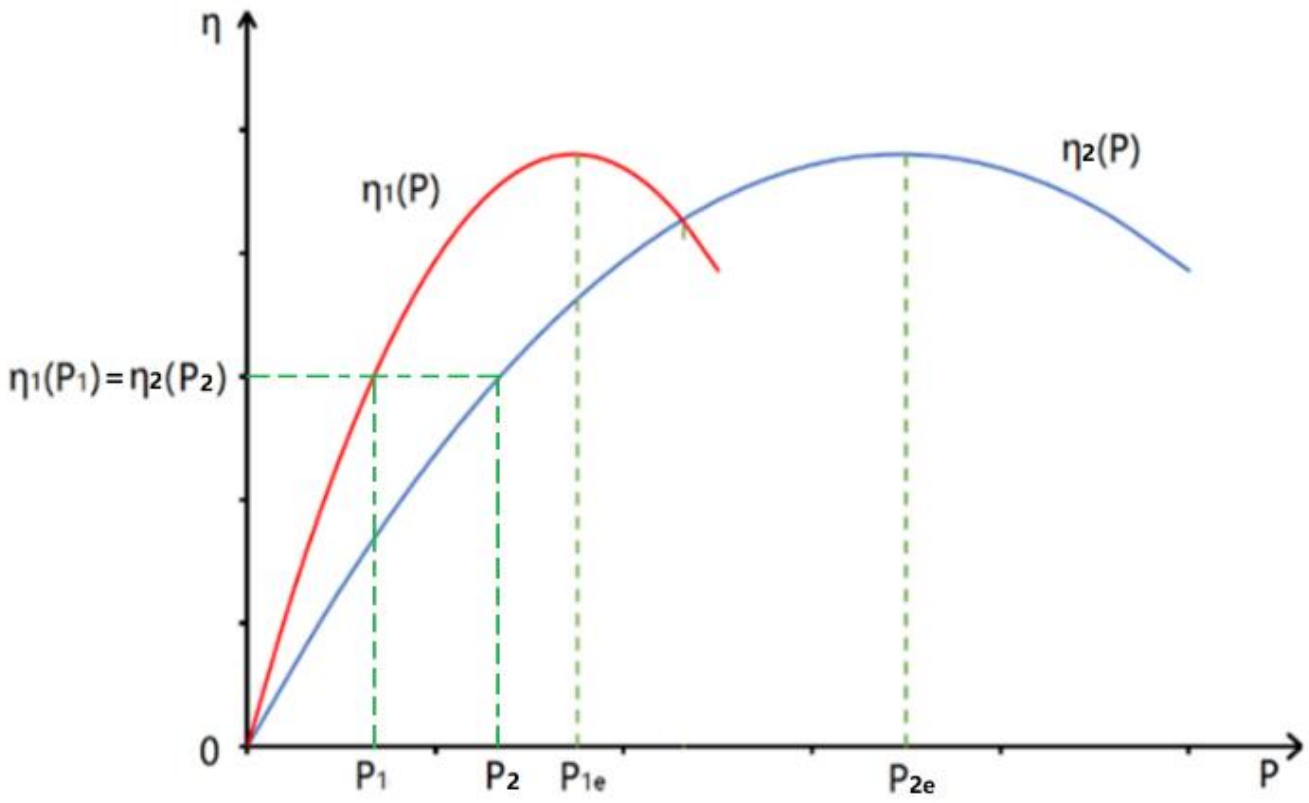

Fig. 3 Energy efficiency curves$\eta_1$(P) and $\eta_2$(P)

b) If $\eta_1(P_1) \neq \eta_2(P_2)$, equation (11) holds, there must be

$$\eta_1(P_1) + P_1\eta_1'(P_1) = \eta_2(P_2) + P_2\eta_2'(P_2) \tag{13}$$

If $\eta_1(P_1) < \eta_2(P_2)$, then $P_1\eta_1'(P_1) > P_2\eta_2'(P_2)$;

If $\eta_1(P_1) > \eta_2(P_2)$, then $P_1\eta_1'(P_1) < P_2\eta_2'(P_2)$;

If $\eta_1'(P_1) > 0$ and $\eta_2'(P_2) < 0$, then $\eta_1(P_1) < \eta_2(P_2)$;

If $\eta_1'(P_1) < 0$ and $\eta_2'(P_2) > 0$, then $\eta_1(P_1) > \eta_2(P_2)$.

For the devices with similar efficiency, it is easy to see

$$\begin{gathered} \eta_2(P_2) = \eta_1(P_1) \\ P_2 = P_1\frac{\eta_1'(P_1)}{\eta_2'(P_2)} \\ P_t = P_1 + P_2 \end{gathered} \tag{14}$$

is an optimization point.

Since $W_t$ is a combination of the efficiency functions of each device, its second derivative has the same sign as the efficiency function of each device, and the second derivative is also less than zero.

$$W_t''(P_1) < 0 \tag{15}$$

$W_t$ has a unique maximum value, and the overall energy efficiency has a unique maximum value also. We have

$$\begin{gathered} \max\eta_{t2}(P_t) = \eta_1(P_1) \\ maxW_t = P_t\eta_1(P_1) \end{gathered} \tag{16}$$

2) n=3

The system has three variables $P_1$, $P_2$ and $P_3$

$$\begin{gathered} P_1 + P_2 + P_3 = P_t \\ P_1 > 0 \\ P_2 > 0 \\ P_3 > 0 \end{gathered} \tag{17}$$

$W_t$ expression is

$$W_t = P_1\eta_1(P_1) + P_2\eta_2(P_2) + P_3\eta_3(P_3) \tag{18}$$

Assuming that $P_3$ is the optimal point and fixed, $P_1$ and $P_2$ are variables, then we have

$$P_1 + P_2 = P_t - P_3 \tag{19}$$

According to the conclusion of n=2, we have

$$\begin{gathered} \eta_2(P_2) = \eta_1(P_1) \\ P_2 = P_1\frac{\eta_1'(P_1)}{\eta_2'(P_2)} \end{gathered} \tag{20}$$

In the same way, assuming that P2 is the optimal point and is fixed, then we have

$$\begin{gathered} \eta_3(P_3) = \eta_1(P_1) \\ P_3 = P_1\frac{\eta_1'(P_1)}{\eta_3'(P_3)} \end{gathered} \tag{21}$$

Assuming that $P_1$ is the optimal point and it has been fixed, we have

$$\begin{gathered} \eta_3(P_3) = \eta_2(P_2) \\ P_3 = P_2\frac{\eta_2'(P_2)}{\eta_3'(P_3)} \end{gathered} \tag{22}$$

It is easy to see

$$\begin{gathered} \eta_1(P_1) = \eta_2(P_2) = \eta_3(P_3) \\ P_2 = P_1\frac{\eta_1'(P_1)}{\eta_2'(P_2)} \\ P_3 = P_1\frac{\eta_1'(P_1)}{\eta_3'(P_3)} \\ P_t = P_1 + P_2 + P_3 \end{gathered} \tag{23}$$

is an optimization point.

Similarly, the second derivative is also less than zero.

$$W_t''(P_1) < 0 \tag{24}$$

$W_t$ has the only maximum value, and the overall efficiency has the only maximum value also. There is

$$\begin{gathered} \max\eta_{t3}(P_t) = \eta_1(P_1) \\ maxW_t = P_t\eta_1(P_1) \end{gathered} \tag{25}$$

3) n=k

Based on the conclusion of n=3, we have

$$\begin{gathered} \eta_1(P_1) = \eta_2(P_2) = \cdots = \eta_k(P_k) \\ P_2 = P_1\frac{\eta_1'(P_1)}{\eta_2'(P_2)} \\ P_3 = P_1\frac{\eta_1'(P_1)}{\eta_3'(P_3)} \\ \cdots \\ P_k = P_1\frac{\eta_1'(P_1)}{\eta_k'(P_k)} \\ P_t = P_1 + P_2 + \cdots + P_k \end{gathered} \tag{26}$$

Similarly, the second derivative of $W_t$ is less than zero.

$$W_t''(P_1) < 0 \tag{27}$$

$W_t$ has the only maximum value, and the overall efficiency has the only maximum value also. We have

$$\begin{gathered} \max\eta_{tk}(P_t) = \eta_1(P_1) \\ maxW_t = P_t\eta_1(P_1) \end{gathered} \tag{28}$$

**Optimal Load Distribution Theorem (Yao Theorem 1)**: For the devices with similar efficiencies, the optimal load distribution method is to keep the operating efficiency of each operating device equal.

$$\eta_1(P_1) = \eta_2(P_2) = \cdots = \eta_n(P_n) \tag{29}$$

**1.5 Optimal Switching Theorem (Yao Theorem 2)**

The optimal methods of load distribution obtained above are all obtained under the assumption that n is already optimal, but is the n optimal?

The total value of input energy A of system C is $P_t$, and there are n operating devices with the different efficiency.

For equipment with similar efficiency, if the number of operating units n is optimal, it must satisfy

$$\eta_1(P_{1,n}) = max\{\eta_1(P_{1,1}), \eta_1(P_{1,2}), \ldots, \eta_1(P_{1,k1}), \eta_1(P_{1,n}), \eta_1(P_{1,k2}), \ldots\} \tag{30}$$

where $P_{1,n}$ is the energy input of the first device when the number of operating devices is n, $P_{1,k1}$ is the energy input of the first device when the number of operating devices is k1, the rest are similar. k1 and k2 are any combination other than the optimal combination of n units this time, and also include other combinations of n units.

Point $P_{1,k1}$ is the point closest to $P_{1,n}$ to the left of point $P_{1,n}$. Point $P_{1,k2}$ is the point closest to $P_{1,n}$ to the right of point $P_{1,n}$, as shown in Fig. 4.

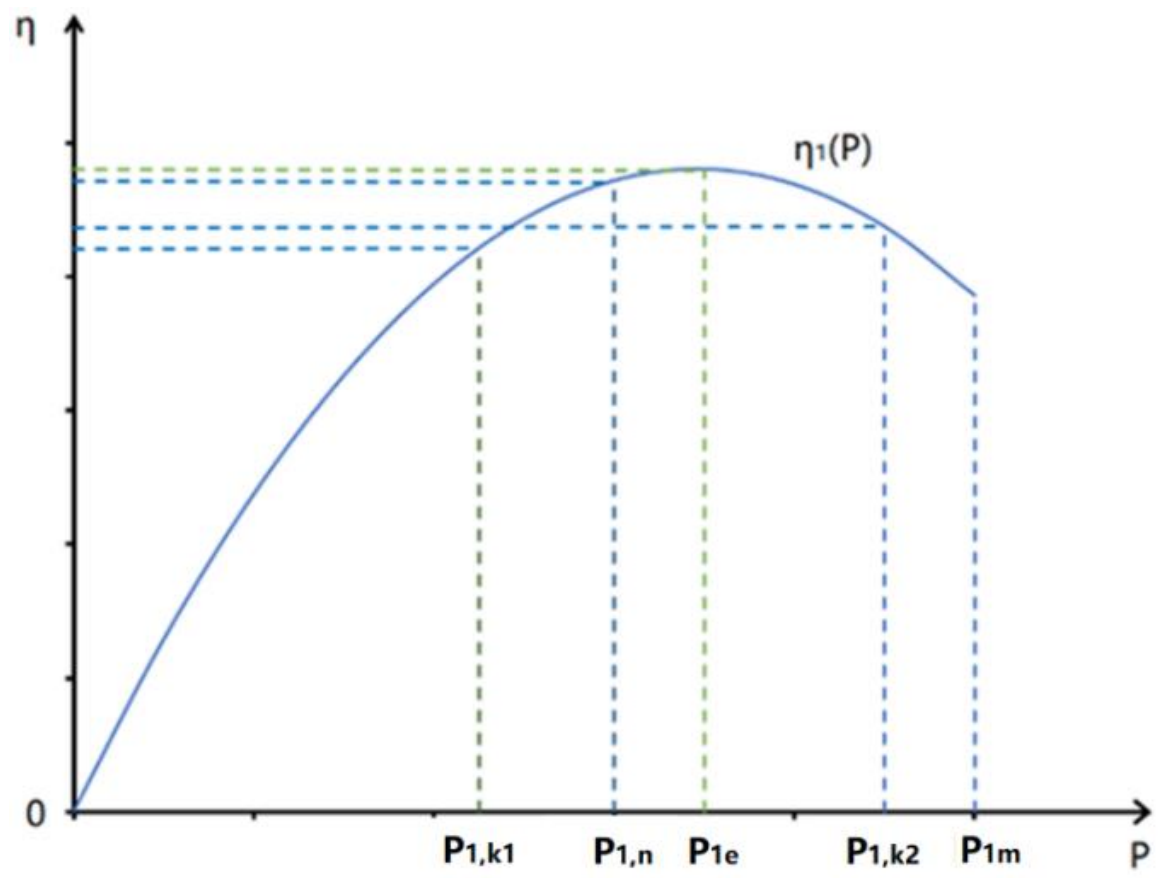

Fig. 4 Energy efficiency comparison curve

$P_t$ increases, $P_{1,n}$ and $P_{1,k1}$ also increase. After point $P_{1e}$, $\eta_1(P_{1,n})$ decreases, but $\eta_1(P_{1,k1})$ increases. When the condition is met

$$\eta_1(P_{1,n}) = \eta_1(P_{1,k1}) \tag{31}$$

the switching point has been reached. If $P_t$ continues to increase, it should switch from n operating devices to $k_1$ operating devices.

If there is no point $P_{1,k1}$ that meets the conditions, the switching point is the $P_{1,n}=P_{1m}$.

$P_t$ decreases, $P_{1,n}$ and $P_{1,k2}$ also decrease. $\eta_1(P_{1,n})$ decreases, but $\eta1(P_{1,k2})$ increases. When the condition is met

$$\eta_1(P_{1,n}) = \eta_1(P_{1,k2}) \tag{32}$$

the switching point has been reached. If $P_t$ continues to decrease, it should switch from n operating devices to $k_2$ operating devices.

**Optimal Switching Theorem (Yao Theorem 2)**: The optimal switching point for the number of operating units is at the point of the equal efficiency or at the maximum output point of the devices.

$$\eta_1(P_{1,n}) = \eta_1(P_{1,k})\ or\ P_{1,n} = P_{1m}\ or\ P_{1,k} = P_{1m} \tag{33}$$

**1.6 Simulation Results**

1. Assume that system C has 2 devices with similar efficiencies.

The efficiency function of device No. 1:

$$\eta_1(P_1) = 0.022P_1 - 0.0001375P_1^2 \tag{34}$$

The efficiency function of device No. 2:

$$\eta_2(P_2) = 0.01647P_2 - 0.00006111P_2^2 \tag{35}$$

The output capacity of device No. 2 is greater than the output capability of device No. 1.

The total output $W_t$ of system C is maximized as

$$\begin{gathered} W_t = \max \textstyle\sum_{i=1}^{2} P_i \eta_i(P_i) \\ s.t.\ \textstyle\sum_{i=1}^{2} P_i = P_t > 0 \end{gathered} \tag{36}$$

In this system, n=1 has two combinations and n=2 has one combination.

When n=1, there is a switching point between the first small device No. 1 and the second large device No. 2. According to the equivalent switching theorem, we have

$$\eta_1(P_t) = \eta_2(P_t) \tag{37}$$

Obtain the switching point $P_{11}$ of small device No. 1 and large device No. 2

$$P_{11} = 96 \tag{38}$$

When n=2, the optimal control method is to maintain

$$\begin{gathered} \eta_1(P_1) = \eta_2(P_2) \\ P_t = P_1 + P_2 \\ P_1\eta_1'(P_1) = P_2\eta_2'(P_2) \end{gathered} \tag{39}$$

Solve for P1 as follows

$$\begin{gathered} P_1 = \frac{2}{5}P_t \\ P_2 = \frac{3}{5}P_t \end{gathered} \tag{40}$$

Switching point $P_{12}$ of No. 2 device and 2 devices satisfies

$$\eta_2(P_t) = \eta_1(\frac{2}{5}P_t) \tag{41}$$

from which we have

$$P_{12} = 150 \tag{42}$$

The maximum total output energy $W_t$ of system C is

$$\max W_t = P_t\eta_1(\frac{2}{5}P_t) \tag{43}$$

The maximum value of the overall operating efficiency $\eta_{t2}(P_t)$ of system C is

$$\eta_{t2}(P_t) = \eta_1(\frac{2}{5}P_t) \tag{44}$$

When $P_{11}<P_t<=P_{12}$ and $P_{12}<P_t$, $W_t''<0$, so the $W_t$ is the maximum value.

The overall operating efficiency is shown in Fig. 5.

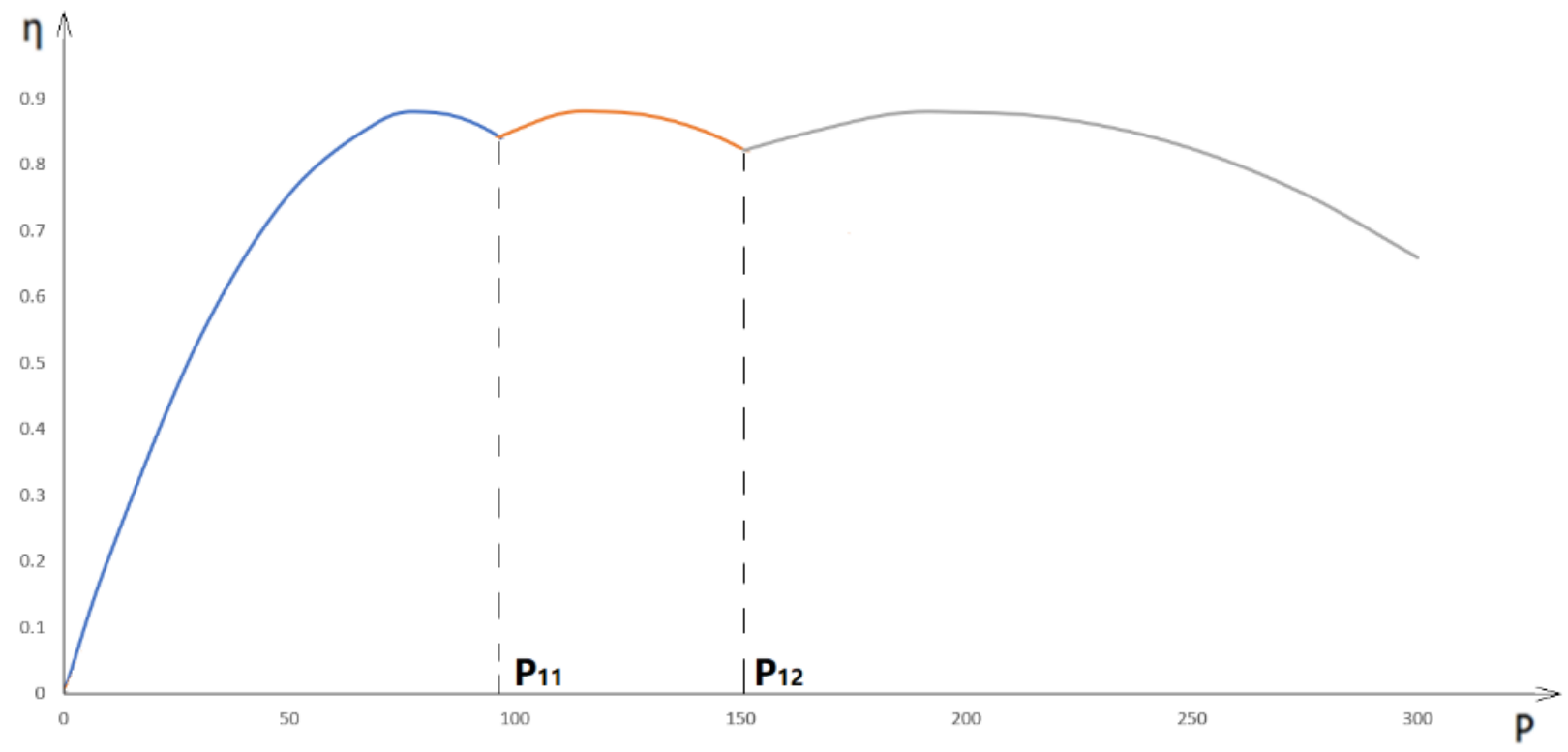

Fig. 5 The overall operating efficiency

2. Assume that system C has 2 devices with no similar efficiencies.

The efficiency function of device No. 1:

$$\eta_1(P_1) = 0.022P_1 - 0.0001375P_1^2 \tag{45}$$

The efficiency function of device No. 2:

$$\eta_2(P_2) = 0.0287P_2 - 0.000233333P_2^2 \tag{46}$$

The output capacity of device No. 2 is smaller than that of device No. 1.

The total output $W_t$ maximization expression of system C is

$$W_t = \max \sum_{i=1}^{2} P_i \eta_i(P_i) \quad s.t. \ \sum_{i=1}^{2} P_i = P_t > 0 \tag{47}$$

In this system, n=1 has two combinations, and n=2 has one combination.

When n=1, there is a switching point between the device No. 1 and the device No. 2. We have

$$\eta_1(P_t) = \eta_2(P_t) \tag{48}$$

Obtain the switching point $P_{11}$ of No. 1 device and No. 2 device

$$P_{11} = 62.61 \tag{49}$$

When $P_t$ is smaller than $P_{11}$, No. 2 device is used. When $P_t$ is greater than $P_{11}$, No. 1 device is used.

The $P_{12}$ is obtained by comparing the values of $P_t\eta_1(P_t)$ and $\max(P_1\eta_1(P_1)+P_2\eta_2(P_2))$.

When $P_t$ is less than $P_{12}$, No.1 device is used

$$P_t\eta_1(P_t) > \max(P_1\eta_1(P_1) + P_2\eta_1(P_2)) \tag{50}$$

When $P_t$ is greater than $P_{12}$, No.1 device and No. 2 device are used

$$P_t\eta_1(P_t) < \max(P_1\eta_1(P_1) + P_2\eta_1(P_2)) \tag{51}$$

$P_1$ and $P_2$ satisfy

$$\eta_1(P_1) + P_1\eta_1'(P_1) = \eta_2(P_2) + P_2\eta_2'(P_2) \tag{52}$$

At the switching point, the efficiency of both is equal

$$\eta_1(P_t) = \frac{P_1\eta_1(P_1)+P_2\eta_2(P_2)}{P_t} = \eta_t(P_t), \quad P_t = P_1 + P_2 \tag{53}$$

Get switching point $P_{12}$

$$P_{12} = 103.36, \quad P_1 = 66.45, \quad P_2 = 36.91 \tag{54}$$

and

$$\eta_1(P_1) = 0.8547$$
$$\eta_2(P_2) = 0.7155 \quad (55)$$
$$\eta_t(P_{12}) = 0.805$$

The overall operating efficiency is shown in Fig. 6.

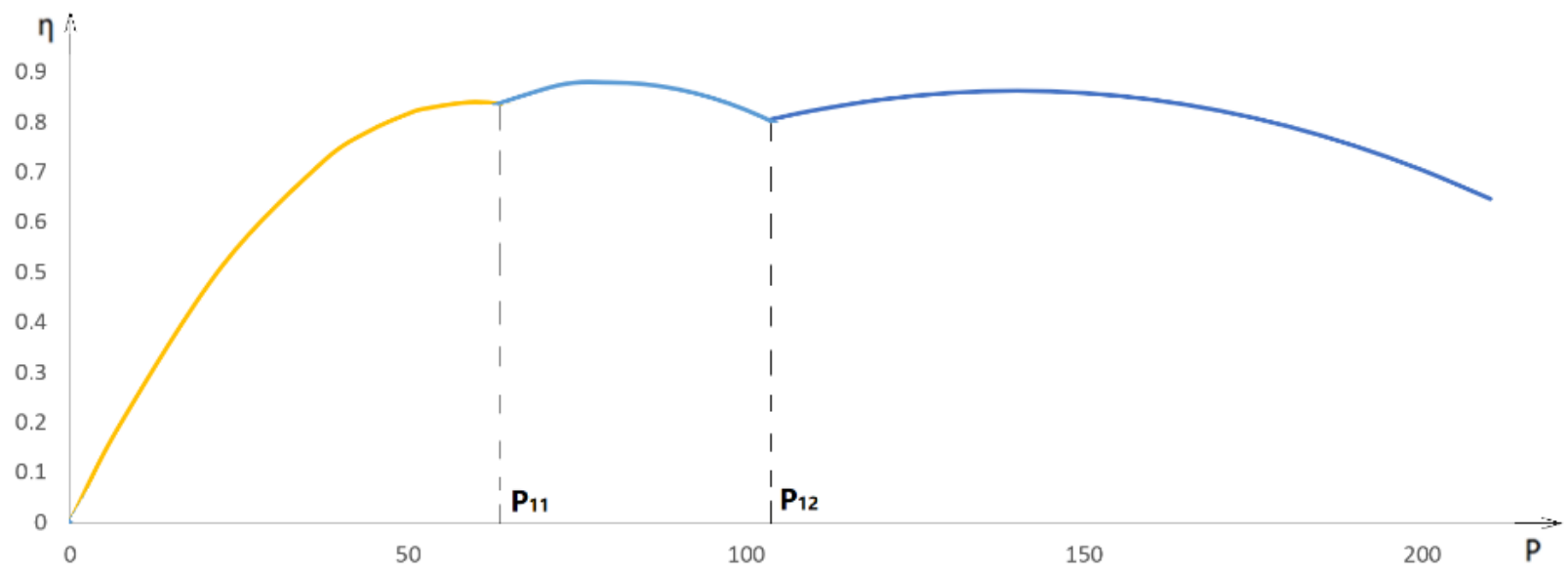

Fig. 6 The overall operating efficiency

### 1.7 Conclusion

There is no need to build an exact mathematical model of the system. According to the characteristics of the efficiency function, a constrained, nonlinear, integer-real hybrid energy efficiency optimization method for multi-unit systems with different efficiencies is proposed.

The optimization method includes two theorems: the optimal load distribution theorem and the optimal switching theorem.

***Optimal load distribution theorem***: The optimal load distribution method for multi-unit systems with different energy efficiencies is to keep the operating efficiency of each operating equipment equal, which is ***Yao's theorem 1***.

$$\eta_1(P_1) = \eta_2(P_2) = \cdots = \eta_n(P_n) \quad (45)$$

***Optimal switching theorem:*** The optimal switching point for the number of operating units is at the equivalent efficiency point or maximum output point of the device, ***Yao Theorem 2.***

$$\eta_1(P_{1,n}) = \eta_1(P_{1,k}) \; or \; P_{1,n} = P_{1m} \; or \; P_{1,k} = P_{1m} \quad (46)$$

We call this optimization method ***quantum optimization method of multi-unit system.***

We call this theory ***energy efficiency predictive theory of multi-unit system.***

These methods have the following advantages:

1) Easy to use, no need to establish an accurate mathematical model of the system;

2) Strong versatility, including linear systems, nonlinear systems, multi-variable systems, time-invariant systems, and time-varying systems;

3) Integer optimization and real number optimization are solved together.

### 1.8 Acknowledgment

In order to solve the problem of efficiency optimization of multi-unit systems, the authors conducted long-term research. We would like to thank Yanfang Zhang from Hebei Automation Company and Bosheng Yao from Beijing IAO Technology Development Company for their valuable help and suggestions during the development and experiment of this theory. I would like to thank my doctoral supervisor Professor Hexu Sun, postdoctoral supervisor Dr. Chengyu Cao, and my supervisors Dr. Qingbin Gao and Dr. Robert X. Gao for their valuable help and advice.